\documentclass[smallextended,envcountsect]{svjour3} 
\smartqed 
\usepackage{amsmath,amssymb}
\usepackage{array}
\usepackage{booktabs}
\usepackage[font={small}]{caption}
\usepackage{cite}
\usepackage{color}
\usepackage{enumitem}
\usepackage{float}
\usepackage{graphicx}
\usepackage{hyperref}
\usepackage{listings}
\usepackage{pbox}
\usepackage{pdflscape}
\usepackage{subfig}
\usepackage{tikz}
  \usetikzlibrary{
    3d,
    arrows, 
    backgrounds,
    positioning, 
    shapes.geometric
  }
  \tikzstyle{special-box} = [draw=black, fill={rgb:red,34;green,139;yellow,34}, 
      opacity=.6, text opacity=1,
      rectangle, rounded corners, inner sep=10pt, inner ysep=10pt]
  \tikzstyle{dsm-box} = [draw=black, fill={rgb:red,140;green,0;yellow,26}, 
      opacity=.6, text opacity=1,
      rectangle, rounded corners, inner sep=10pt, inner ysep=10pt]
  \tikzstyle{test-box} = [draw=black, fill={rgb:red,70;green,130;blue,180},
      opacity=.6, text opacity=1,
      diamond, aspect=3, inner sep=1pt, inner ysep=1pt]
  \tikzset{
      descr/.style={
          fill=white,
          inner sep=2.5pt
      }
  }

\graphicspath{{figures/}}
\floatstyle{boxed}
\newfloat{procedure}{thp}{lop}

\newtheorem{myAlgorithm}{Procedure}

\usepackage{xspace}

\newcommand{\CalC}[1]{\ensuremath{\mathcal{C}^{#1}}\xspace}

\newcommand{\Man}{\ensuremath{\mathcal{M}}\xspace}
\newcommand{\MyExp}[2]{\ensuremath{\mathbf{Exp}_{#1}(#2)}\xspace}
\newcommand{\GLn}[1]{\ensuremath{\mathtt{GL}(#1)}\xspace}
\newcommand{\GLnp}{\ensuremath{\mathtt{GL}^{+}(n)}\xspace}
\newcommand{\Gnk}[2]{\ensuremath{\mathtt{G}(#1,#2)}\xspace}
\newcommand{\la}[1]{\ensuremath{\mathfrak{#1}}\xspace}
\newcommand{\on}[1]{\ensuremath{\mathfrak{o}(#1)}\xspace}
\newcommand{\On}[1]{\ensuremath{\mathtt{O}(#1)}\xspace}

\newcommand{\Rn}[1]{\ensuremath{\mathbb{R}^{#1}}\xspace}
\newcommand{\SEn}[1]{\ensuremath{\mathtt{SE}(#1)}\xspace}
\newcommand{\sen}[1]{\ensuremath{\mathfrak{se}(#1)}\xspace}
\newcommand{\SLn}{\ensuremath{\mathtt{SL}(n)}\xspace}
\newcommand{\sln}{\ensuremath{\mathfrak{sl}(n)}\xspace}
\newcommand{\Sn}[1]{\ensuremath{\mathcal{S}^{#1}}\xspace}
\newcommand{\SOn}[1]{\ensuremath{\mathtt{SO}(#1)}\xspace}
\newcommand{\son}[1]{\ensuremath{\mathfrak{so}(#1)}\xspace}
\newcommand{\SPDn}{\ensuremath{\mathtt{SPD}(n)}\xspace}
\newcommand{\spdn}{\ensuremath{\mathfrak{spd}(n)}\space}
\newcommand{\TMan}[1]{\ensuremath{\mathcal{T}_{#1}\mathcal{M}}\xspace}
\newcommand{\Tn}[1]{\ensuremath{\mathtt{T}(#1)}\xspace}
\newcommand{\tn}[1]{\ensuremath{\mathfrak{t}(#1)}\xspace}
\newcommand{\UPn}{\ensuremath{\mathtt{UP}(n)}\xspace}
\newcommand{\Vnk}[2]{\ensuremath{\mathtt{V}(#1,#2)}\xspace}

\journalname{JOTA}

\begin{document}

\title{Direct Search Methods on Reductive Homogeneous Spaces}

\author{David W. Dreisigmeyer
\thanks{Communicated by Alexandru Krist\'{a}ly.}
}

\institute{
  David W. Dreisigmeyer,  Corresponding author \at
  United States Census Bureau \\
  Center for Economic Studies \\
  Suitland, MD \\
  Department of Electrical and Computer Engineering \\
  Colorado State University \\
  Fort Collins, CO \\
  david.wayne.dreisigmeyer@census.gov
}

\date{Received: date / Accepted: date}

\maketitle

\begin{abstract}
	Direct search methods are mainly designed for use in problems with no equality constraints.  
    However there are many instances where the feasible set is of measure zero in the ambient space and no mesh point lies within it.  
    There are methods for working with feasible sets that are (Riemannian) manifolds, but not all manifolds are created equal.  
    In particular, reductive homogeneous spaces seem to be the most general space that can be conveniently optimized over.  
    The reason is that a `law of motion' over the feasible region is also given.  
    Examples include \Rn{n} and its linear subspaces, Lie groups and coset manifolds such as Grassmannians and Stiefel manifolds.  
    These are important arenas for optimization, for example, in the areas of image processing and data mining.  
    We demonstrate optimization procedures over general reductive homogeneous spaces utilizing maps from the tangent space to the manifold.  
    A concrete implementation of the probabilistic descent direct search method is shown.  
    This is then extended to a procedure that works solely with the manifold elements, eliminating the need for the use of the tangent space.
\end{abstract}

\section{\label{sec:introduction}Introduction}
Direct search methods are typically designed to work in a full-dimensional subset of \Rn{n}. 
So, if we only have (at most) inequality constraints present, these methods are potentially useful procedures to employ, e.g., \cite{art:audet-2004, art:audet-2006, art:coope-2000, art:dennis-2004, tech:finkel-2004, art:price-2003}.
This is especially true for black-box optimization or cases where derivative information is not available.  
However the presence of even one nonlinear equality constraint defining a feasible set of measure zero creates difficulties.  
In this case the probability of a mesh point being feasible is zero.  
It is shown how this difficulty can be avoided in some particularly convenient settings.

The mesh adaptive direct search (MADS) algorithms \cite{art:abramson-2009, art:audet-2006} provide typical illustrations.  
The probability that a point on the mesh used in MADS also lies in the feasible set is zero.  
So every point that MADS considers will likely be infeasible.  
Filter \cite{art:audet-2004, art:fletcher-2002} and progressive barrier \cite{art:audet-2009} techniques may be able to alleviate this situation if one is willing to violate the equality constraints somewhat.  
Also, the augmented Lagrangian pattern search algorithm in \cite{art:lewis-2002} extends the method in \cite{art:conn-1991} to non-smooth problems with equality constraints.  
We have in mind something quite different from these methods.  

A recent example of the type of problem we will be addressing is \cite{art:audet-2015}.  
The method presented there handles linear equalities in the MADS algorithm \cite{art:audet-2006}.  
The method presented here will deal with subsets of \Rn{n} that can have a much more complicated structure.

Optimization is naturally carried out in \Rn{n}, over Lie groups \cite{book:Stillwell-2008} and within other geometric objects like Grassmannians \cite{art:edelman-1999}.  
A common theme running through these is that of reductive homogeneous (often symmetric) spaces \cite{book:Sharpe-1997}.  
Globally these have a simple structure and convenient ways to move over them.  
Exploring how to optimize within them has been a worthwhile pursuit \cite{art:edelman-1999, art:Smith-1994}. 

Following Sharpe \cite{book:Sharpe-1997}, Euclidean geometry was generalized through two routes: Riemannian geometry and Klein geometries.  
Euclidean geometry has: 1) a space \Rn{n}, 2) with a symmetry Lie group, the special Euclidean group, that allows us to reach any point in \Rn{n} from any other point, and 3) a metric.  
Riemannian geometry generalizes Euclidean geometry by dropping the symmetry Lie group condition allowing arbitrary curvature of our underlying space.  
Locally the manifold looks like \Rn{n}.  
However the symmetry group provides a convenient way to navigate within a manifold.  
A Klein geometry generalizes the symmetry groups that act on the underlying spaces: it specifies a space and a way to move in that space.

Generalizing optimization over \Rn{n} has tended to focus on manifolds with a Riemannian structure \cite{book:absil-2008}.  
Importantly, Sard's theorem tells us that equality constraints can typically be taken to define a Riemannian manifold.  
The problem here is that the global structure of this type of manifold can be complicated in which case just moving over them without wandering off of their surface may be difficult.

We will extend direct search methods to constrained problems when the constraints are specified by a Klein geometry.  
The Klein geometry will implicitly define a subset of the ambient space \Rn{n} of our optimization variables.   
This subset is the feasible set, possibly further constrained with inequality constraints.  
The advantage of using a Klein geometry is that it additionally gives us information on how to move over the feasible set.

The current paper is laid out as follows.  
Section~\ref{sec:manifolds} provides a brief overview of the necessary differential geometry.  
Lie groups and Lie algebras are introduced in Section~\ref{sec:lie-first}.  
Section~\ref{sec:background} illustrates our methodology through a series of examples.  
We return for a look at Lie groups and algebras in Section~\ref{sec:lie-second} before presenting the optimization procedure in Section~\ref{sec:convergence}.  
In Section~\ref{sec:in-group} we investigate an optimization routine where one can remain strictly in the Lie group and avoid any use of the Lie algebras.  
This is perhaps the highlight of the current paper and represents an alternate way of doing optimization.  
Some data mining applications of the new procedures are looked at in Section~\ref{sec:applications}.  
A discussion follows in Section~\ref{sec:discussion}.  
Not all possibilities are dealt with: the focus is on what would be considered the easiest spaces that would most likely be used in practice.

\section{\label{sec:manifolds}Manifolds}
For us a manifold \Man can be thought of as a subset of \Rn{n} that looks locally like \Rn{m}, $n > m$.  
To formalize this, let $\mathcal{I}$ be an index set and for each $i \in \mathcal{I}$ let $\mathcal{V}_{i} \subset \Man$ be an open subset where $\bigcup_{i \in \mathcal{I}} \mathcal{V}_{i} = \Man$.  
Additionally, let there be homeomorphisms $\phi_{i} : \mathcal{V}_{i} \rightarrow \Rn{m}$.  
Each pair $(\mathcal{V}_{i}, \phi_{i})$ is called a \emph{chart} and the set $\{ (\mathcal{V}_{i}, \phi_{i}) \}_{i \in \mathcal{I}}$ is called an \emph{atlas}.  
We now have local coordinates etched onto \Man.

When two open subsets $\mathcal{V}_{i}$ and $\mathcal{V}_{j}$ of \Man overlap there will be a homeomorphic transition function $\tau_{ij} : \Rn{m} \rightarrow \Rn{m}$ on the overlap $\mathcal{V}_{i} \cap \mathcal{V}_{j}$ defined by $\tau_{ij} = \phi_{j} \circ \phi_{i}^{-1}$.  
The transition functions allow us to stitch together the charts into a coherent whole.  
This then defines a \emph{topological manifold}.

Extra structure can be placed on the manifold by placing additional requirements on \Man.  
A typical requirement is that the transition functions are differentiable.  
In this case we have a \emph{differentiable manifold}.  
If the $\tau_{ij}$ are all \CalC{k}, $k \geq 1$, then \Man is said to be a \CalC{k}-manifold.  
If \Man is a \CalC{\infty}-manifold then it is called \emph{smooth}.  
Every \CalC{k}-manifold can be made smooth \cite{book:hirsch-1997}.

At every point $p \in \Man \subset \Rn{n}$ a tangent space \TMan{p} can be attached.  
If an inner product $h: \TMan{p} \times \TMan{p} \rightarrow \mathbb{R}$ is defined for every $p \in \Man$ and varies smoothly over \Man then we have a \emph{Riemannian manifold}.  
The function $h(\cdot,\cdot)$ is called a \emph{metric} and with it notions of lengths and angles in \TMan{p} are available over \Man.  
We will need surprisingly little direct use of any Riemannian structure.

Manifolds often enter an optimization problem through equality constraints.  
\begin{definition}[Regular level sets]
\label{def:level-sets}
Let $\mathbf{g}: \Rn{n + m} \rightarrow \Rn{m}$.  
For the set $\Man = \{ \mathbf{x}:\mathbf{g}(\mathbf{x}) = \mathbf{c} \}$ assume $\nabla \mathbf{g}(\mathbf{x})$ is full rank.  
Then \Man is a \emph{regular level set} of $\mathbf{g}(\mathbf{x})$.  
Additionally, the null space of $\nabla \mathbf{g}(\mathbf{x})$ coincides with the tangent space to \Man.
\end{definition}
If the function defining the equality constraints is sufficiently differentiable then by Sard's theorem almost every level set is a (Riemannian) manifold.
\begin{theorem}[Regular level set manifolds and Sard's theorem \cite{book:helgason}]
\label{theo:sard}
Let $\mathbf{g}: \Rn{n + m} \rightarrow \Rn{m}$ be a \CalC{2} function.  
Then every regular level set of $\mathbf{g}(\mathbf{x})$ is an $n$-dimensional manifold.  
Further, if $\mathbf{g}(\mathbf{x})$ is a \CalC{n + 1} function then almost every level set is a regular level set.
A regular level set acquires a metric from its embedding in \Rn{n + m} by restricting the standard Euclidean metric to the tangent space \TMan{p} given by $\mathrm{null}(\nabla \mathbf{g}(\mathbf{p}))$ for all $p \in \Man$.
\end{theorem}

On a Riemannian manifold \Man a \emph{geodesic} is the shortest path between two (nearby) points on \Man and generalizes the notion of a straight line.  
The exact way this is calculated for a Riemannian manifold is not important for us.  
What is important though is when any arbitrary point $a \in \Man$ can be reached from every other point $b \in \Man$ by a geodesic.  
In such a case the manifold \Man is said to be \emph{geodesically complete}.  
This has many nice properties which are covered in more detail in \cite{misc:dreisigmeyer-2006}.  
Not all of the manifolds considered herein will be geodesically complete.  
When they are not considerably more effort may be required to implement a practical optimization routine.

\section{\label{sec:lie-first}Lie Groups and Lie Algebras}
A \emph{Lie group} $G$ is a smooth manifold with the additional structure of being a group.  
For $g,h \in G$ the function $\phi(g,h) = g^{-1}h$ must additionally be smooth.  
A \emph{matrix Lie group} is defined as any subgroup of the group of general linear transformations \GLn{n} of \Rn{n}, the Lie group of invertible $n$-by-$n$ matrices.  
We will always work with matrix Lie groups and simply say Lie group instead of matrix Lie group.

A major advantage of having the group structure is the ability to easily move over the manifold.  
Every element $g \in G$ defines a diffeomorphism of $G$ by the operation of \emph{left translation} defined by $L_{g}(h) = g \odot h$ for every $h \in G$.  
Left translation has moved the element $h$ to the new point $g \odot h$ on the manifold.
Often $g \odot h$ is simply matrix multiplication but this is not always the case.
For instance symmetric positive definite matrices would have $g \odot h = g^{1/2} h g^{1/2}$.
We will drop the $\odot$ with the understanding that the correct multiplication is used when it differs from the usual matrix multiplication.

A \emph{path} in $G$ is given by $L_{h_{1}} \circ L_{h_{2}} \circ \cdots \circ L_{h_{n}} (e) = h_{1} \cdots h_{n} = L_{h_{1} \cdots h_{n}}(e)$.  
The \emph{path-connected identity component} of a Lie group $G$ is the subgroup of those elements of the full group that can be connected to the identity element by a path.  
We always assume that this is the component that is being used unless stated otherwise.

Since a Lie group is a manifold it will have a tangent space \TMan{e} at the identity.  
\TMan{e} is called the \emph{Lie algebra} \la{g} of the Lie group $G$.  
A Lie algebra is a vector space and is a much easier object to work with in a typical optimization procedure.  
For example, we can lay out a mesh for MADS in a Lie algebra.  
For a matrix group the elements of the Lie algebra are also given by matrices.  
The \emph{Lie bracket} on the Lie algebra is defined as $[\omega_{1}, \omega_{2}] = \omega_{1}\omega_{2} - \omega_{2}\omega_{1}$ for $\omega_{1}, \omega_{2} \in \la{g}$.

There is a mapping $\exp: \la{g} \rightarrow G$ called the \emph{exponential mapping}.  
For matrix Lie groups this corresponds to the matrix exponential.  
The following is a key fact for Lie groups:
\begin{theorem}[Path from the identity \cite{book:Stillwell-2008}]
\label{th:lie-path}
Every element of a path-connected matrix Lie group is given by a finite product of exponents of elements in the corresponding Lie algebra.  
That is, if $g \in G$ then $g = \exp[\omega_{1}] \cdots \exp[\omega_{n}]$ for some $\omega_{i} \in \la{g}$.
\end{theorem}
The exponential map $\exp[t \omega]$, $t \in \mathbb{R}$ and $\omega \in \la{g}$, is a geodesic that sometimes, but not always, corresponds to the Riemannian geodesic.  
When it doesn't correspond we will use the Lie group exponential.  
An important special case is when the Lie group $G$ is compact (as a manifold).  
Then the exponential map is also a Riemannian geodesic and the manifold is geodesically complete.

Given a Lie group $G$ and a subgroup $H \subset G$ we can form the \emph{coset} $G/H$ where every element in $G/H$ is an equivalency class $gH = \{gh:h \in H\}$ for every $g \in G$.  
We identify $g_{1}, g_{2} \in G$ as the same element if $g_{1}H = g_{2}H$ as sets.  
The space $G/H$ is also a smooth manifold, called the \emph{coset manifold}.

\section{\label{sec:background}Examples}
It is important at this point to lay out what is needed to effectively perform an optimization procedure.  
What is it about \Rn{n} that needs to be retained and what can be thrown away?  
There are three vital things that need to be done:
\begin{enumerate}
	\item move in a specified direction;
	\item a specified distance, and;
	\item evaluate a function at every point.
\end{enumerate}
If we can do these we can optimize.  
If we can do these efficiently we can optimize efficiently.

For Riemannian geometries all of the three requirements are satisfied.  
We can move a specified length along a specified direction on the manifold.  
Having functions defined on a Riemannian manifold also causes no difficulties.  
The difference from \Rn{n} is that the path we travel along now, our 'direction', is not necessarily a straight line.  
Instead it's a geodesic, the straightest possible line we could hope for and still remain on the manifold.  
But this is the difficulty: numerically traveling along the geodesic without straying is hard.

There are some cases of a Riemannian manifold where it is easy to stay on the geodesic.  
\Rn{n} is the canonical example where the geodesics are simply straight lines and the Euclidean metric gives us our sense of distance.  
However, it is possible for other Riemannian manifolds to have closed-form solutions for their geodesics and this makes optimizing over them a relatively straightforward procedure.  

If the surface we are optimizing over has a simple global structure then we can hope that performing the optimization will be easy.  
A particular example of this is when we know the Lie group of symmetries of the manifold.  
In this case the group elements specify how to move over the manifold.  
Below we examine this through a few examples of progressive complexity.  
There will be results and definitions that, while expressed in a specific example, are generic for the spaces and Lie groups we'll consider.

\subsection{\label{subsec:se3}The Special Euclidean Group \SEn{3}}
Our first case will be \Rn{3}.  
Here the orientation preserving symmetry group is given by the special Euclidean group \SEn{3} that includes all of the rotations and translations of \Rn{3}.  
That is, the space looks the same under any transformation in \SEn{3}.  
The typical way we would represent these transformations is by augmenting a vector $\mathbf{x} \in \Rn{3}$ with a final $1$ to have $\hat{\mathbf{x}}^{T} = [\ \mathbf{x}^{T}\ 1\ ]$.  
A matrix representation $M$ of a group element $g \in \SEn{3}$ is then given by
\begin{equation}
\label{eq:sen-matrix-rep}
	M = 
    	\left[
        	\begin{array}{cc}
         		O 				&	\mathbf{v} \\
            	\mathbf{0}^{T}	&	1
          	\end{array}
    	\right] \mbox{,}
\end{equation}
where $O \in \SOn{3}$ is a $3$-by-$3$ orthogonal matrix with $\det(O) = 1$ and $\mathbf{0}, \mathbf{v} \in \Rn{3}$.  
Then an explicit representation of the action of $g \in \SEn{3}$ on $\mathbf{x} \in \Rn{3}$ is given by $M \hat{\mathbf{x}}$.  

With \SEn{3} we can reach any point $\mathbf{y} \in \Rn{3}$ from any other point $\mathbf{x} \in \Rn{3}$.  
However this is not unique since we can always first perform a rotation around $\mathbf{x}$ and then translate from $\mathbf{x}$ to $\mathbf{y}$.  
By `rotation around $\mathbf{x}$' we mean translating $\mathbf{x}$ to the origin, performing a rotation and then translating back to $\mathbf{x}$.  
There's no reason to consider these (infinite number of) different ways of going from $\mathbf{x}$ to $\mathbf{y}$ as being different.  
Instead we form equivalency classes where $M_{1} \sim M_{2}$ when
\begin{equation}
\label{eq:sen-matrix-equiv}
	M_{1} = M_{2}
    	\left[
        	\begin{array}{cc}
         		O 				&	\mathbf{0} \\
            	\mathbf{0}^{T}	&	1
          	\end{array}
    	\right] \mbox{,}
\end{equation}
for $O \in \SOn{3}$.  
Notice that the matrix that multiplies $M_{2}$ on the right in (\ref{eq:sen-matrix-equiv}) is also a representation of $O \in \SOn{3}$ since the `extra' $1$ is irrelevant.  
It is also a member of \SEn{3} where there is no translation.  
We see that $\SOn{3} \subset \SEn{3}$ and additionally that $\SEn{3} = \SOn{3} \ltimes \Tn{3}$, the semi-direct product of rotations and translations in \Rn{3}.  
From this we can form the coset manifold $\SEn{3} / \SOn{3} = \Tn{3} \simeq \Rn{3}$, getting back a copy of what we started with.  
The optimization proceeds in the coset manifold.

A \textit{Klein geometry} $(G,H)$ consists of a Lie group $G$, a subgroup $H \subset G$ and a space (manifold) $\Man = G / H$ over which $G$ acts \emph{transitively}: for every $p,q \in \Man$ there is a $g \in G$ such that $gp = q$.  
Alternatively, a Klein geometry can be thought of as a space and a law for moving over this space.  
For the Klein geometry $(\SEn{3}, \SOn{3})$ the following conditions for the Lie algebras $\sen{3} = \son{3} \oplus \tn{3}$ hold:
\begin{equation}
\label{eq:sen-sym-space}
	\left[ \son{3} , \son{3} \right] \subseteq \son{3} \mbox{, }
    \left[ \son{3} , \tn{3} \right] \subseteq \tn{3} \mbox{ and }
    \left[ \tn{3} , \tn{3} \right] \subseteq \son{3} \mbox{, }
\end{equation}
where $[\cdot,\cdot]$ is the Lie bracket.  
The first condition is that \son{3} is itself a Lie algebra.  
The second condition tells us that the \tn{3} is invariant under the actions of \son{3}.  
The final condition states that the Lie bracket of two infinitesimal translations is a rotation.  
As these conditions are sequentially imposed we move from a \textit{homogeneous space} to a \textit{reductive homogeneous space} to a \textit{symmetric space}.  

We are still missing a sense of moving along a straight line in \Rn{3} by a specified amount.  
One thing that makes reductive homogeneous spaces convenient to work with is that if the Klein geometry is given by $(G,H)$, and $\la{g} = \la{h} \oplus \la{m}$ with $[\la{h},\la{h}] \subseteq \la{h}$ and $[\la{h},\la{m}] \subseteq \la{m}$, then we can ignore the action of $\la{h}$ and use only $\la{m}$ in our optimization.  
What this means is we can use the same exponential map defined on $\la{g}$ and simply restrict it to $\la{m}$ to find the exponential map to $G/H$.  
It is the reductive condition that allows this.

In the case of \tn{3} the matrices in the algebra are of the form
\begin{equation}
\label{eq:tn3-matrix-rep}
	\tau = 
    	\left[
        	\begin{array}{cc}
         		0_{3} 				&	\mathbf{v} \\
            	\mathbf{0}^{T}	&	0
          	\end{array}
    	\right] \mbox{,}
\end{equation}
and $\exp[\tau] = I_{4} + \tau$.  
This explicitly shows that $\Tn{3} \simeq \Rn{3}$ through the identification $\tau \leftrightarrow \mathbf{v}$.  
We can travel along a straight line from from the point $T_{0} \in \Tn{3}$ using $T_{0}\exp[t \tau]$, $t \in \mathbb{R}$.  
As $t$ increases we are traveling further along this geodesic.

A symmetric space is always geodesically complete, which does not necessarily hold for the other homogeneous spaces \cite{book:Sternberg-2012}.  
Optimizing over a geodesically complete space is easier since this allows one to use a global \textit{pullback procedure} from the manifold to the tangent space as covered in \cite{misc:dreisigmeyer-2006}.  
This means we can pullback the objective function $f: \Man \rightarrow \mathbb{R}$, and any inequality constraints, from the manifold \Man to the Lie algebra \la{m} using the function composition $f \circ \exp : \la{m} \rightarrow \mathbb{R}$.  
Since the Lie algebra is a normed vector space and the exponential map is smooth, any convergence results about a direct search method in \Rn{n} automatically hold without modification in the Lie algebra.  
In this case optimizing over the underlying space becomes a standard optimization in Euclidean space apart from a more complicated exponential function.

\subsection{\label{subsec:gnk}The Grassmannian \Gnk{4}{2}}
Our second example will still be a symmetric space but more abstract.  
Here we will look at the Grassmannian \Gnk{4}{2} of all two-dimensional planes in \Rn{4}.  
We start with the Lie group \On{4}, which is the symmetry group of the $3$-sphere $\mathcal{S}^{3} \subset \Rn{4}$.  
A fixed plane in \Rn{4} can be specified by the subspace normal to it.  
The subtlety for \Gnk{4}{2} is that we don't want to differentiate between any orthonormal bases that span the plane or the subspace normal to it.  
Then $\Gnk{4}{2} \simeq \On{4} / (\On{2} \times \On{2})$.  
The first \On{2} ignores any changes in the orthonormal basis for the subspace normal to the plane; the second ignores any changes in the orthonormal basis for the plane.  

For \Gnk{4}{2}, $\on{4} = \la{h} \oplus \la{m}$ where
\begin{eqnarray}
\label{eq:on3-matrix-rep}
	\omega & = & \eta + \mu \nonumber \\
    	& = &
        \left[
        	\begin{array}{cc}
         		A		&	0 \\
            	0		&	B
          	\end{array}
    	\right]
        +
        \left[
        	\begin{array}{cc}
         		0 		&	-C \\
            	C^{T}	&	0
          	\end{array}
    	\right]
        \mbox{,}
\end{eqnarray}
with $A$ and $B$ skew-symmetric, $\eta \in \la{h}$ and $\mu \in \la{m}$.  
Since
\begin{equation}
\label{eq:vnk-sym-space}
	\left[ \la{h} , \la{h} \right] \subseteq \la{h} \mbox{, }
    \left[ \la{h} , \la{m} \right] \subseteq \la{m} \mbox{ and }
    \left[ \la{m} , \la{m} \right] \subseteq \la{h} \mbox{, }
\end{equation}
\Gnk{4}{2} is a symmetric space and hence geodesically complete.  
Then the map $\exp[\mu]$, $\mu \in \la{m}$, allows us to reach any point in \Gnk{4}{2} from any other point.  
As for \Rn{3},  we can travel along a geodesic from a point $G_{0} \in \Gnk{4}{2}$ using $G_{0}\exp[t \mu]$, $t \in \mathbb{R}$.

Here the exponential mapping from the Lie algebra to the manifold is significantly more complex than that for the \Rn{3} case \cite{art:edelman-1999}.  
Table~\ref{tbl:metrics_exps} gives the solution for \Gnk{n}{k} along with some other Lie groups and reductive homogeneous spaces.  
Finding a way to efficiently calculate the exponential map for a reductive homogeneous space is the key to being able to efficiently optimize on it.  
For a Klein geometry $(G,H)$ it is enough to know the exponential for the Lie group $G$ if the space $G/H$ is a reductive homogeneous space.  
For us this exponential map is always the matrix exponential.

\subsection{\label{subsec:vnk} The Stiefel Manifold \Vnk{4}{2}}
Our next example is the Stiefel manifold \Vnk{4}{2} of two-dimensional planes in \Rn{4} where the orthonormal bases used to span a two-dimensional plane are distinguished.  
The orthonormal bases that span the normal space to a plane are not important however.  
So $\Vnk{4}{2} \simeq \On{4} / \On{2}$.  Now $\on{4} = \la{h} \oplus \la{m}$ is given by
\begin{eqnarray}
\label{eq:on3-matrix-rep2}
	\omega & = & \eta + \mu \nonumber \\
    	& = &
        \left[
        	\begin{array}{cc}
         		0		&	0 \\
            	0		&	B
          	\end{array}
    	\right]
        +
        \left[
        	\begin{array}{cc}
         		A 		&	-C \\
            	C^{T}	&	0
          	\end{array}
    	\right]
        \mbox{,}
\end{eqnarray}
where $A$ and $B$ are skew-symmetric, $\eta \in \la{h}$ and $\mu \in \la{m}$.  

This is our first example of a reductive homogeneous space that is not a symmetric space:
\begin{equation}
\label{eq:vnk-red-hom-space}
	\left[ \la{h} , \la{h} \right] \subseteq \la{h} \mbox{ and }
    \left[ \la{h} , \la{m} \right] \subseteq \la{m} \mbox{ but }
    \left[ \la{m} , \la{m} \right] \not\subseteq \la{h} \mbox{. }
\end{equation}
So we can't say that \Vnk{4}{2} is geodesically complete by appealing to it being a symmetric space.  
But it is compact and can be made Riemannian so the Hopf-Rinow theorem tells us it is geodesically complete \cite{book:helgason}.  
The same as for \Rn{3} and \Gnk{4}{2}, a geodesic from a point $V_{0} \in \Vnk{4}{2}$ is given by $V_{0}\exp[t \mu]$, $t \in \mathbb{R}$, which is both a Riemannian geodesic and a matrix exponential map.  
Table~\ref{tbl:metrics_exps} has the closed-form solution for the geodesics.

\subsection{\label{subsec:subspace} A Fixed Subspace}
This is a simple example of how equality constraints can be restated as a Klein geometry.  
(See \cite{art:audet-2015} also for an implementation of this sort of procedure.)  
Suppose we have a linear constraint $A^{T} \mathbf{x} = \mathbf{0}$ where $A^{T} : \Rn{n + m} \rightarrow \Rn{m}$.  
The subspace that $\mathbf{x} \in \Rn{n + m}$ must lie in is given by $\mathrm{null}(A^{T})$ and $A$ gives the normal directions to the manifold.  Let $Q \in \Rn{(n + m) \times n}$ be a full rank matrix such that $A^{T}Q = 0$ so that Q spans the null space.  
The linear subspace represented by $\mathrm{span}(Q)$ defines the manifold $\mathcal{Q}$.

Assume now that we additionally have $A^{T} A = I_{m}$ and $Q^{T} Q = I_{n}$, which can always be done.  
The symmetry group associated with $\mathcal{Q}$ is the subgroup of \SEn{n + m} given by matrices of the form
\begin{equation}
\label{eq:fix-subspace-G}
	M = 
    	\left[ 
        	\begin{array}{cc}
         		O 				&	Q \mathbf{v} \\
            	\mathbf{0}^{T}	&	1
          	\end{array}
    	\right] \mbox{ and }
	O =
    	U
        \left[
          \begin{array}{cc}
          W & 0 \\
          0 & V
          \end{array}
        \right]
        U^{T}
\end{equation}
where $U = [A\ Q] \in \SOn{n + m}$, $W \in \SOn{m}$, $V \in \SOn{n}$ and $\mathbf{v} \in \Rn{n}$.  
It's easy to check that $M : \mathcal{Q} \rightarrow \mathcal{Q}$ and that all such $M$ form a Lie group $G$.  

Now we have a space (manifold) given by $\mathcal{Q}$ and a law for moving over $\mathcal{Q}$ given by $G$.  
This defines a Klein geometry.  
The last piece we would like is the subgroup $H \subset G$ that will \emph{stabilize} $\mathbf{0} \in \mathcal{Q}$: for all $h \in H$ we must have $h \mathbf{0} = \mathbf{0}$.  
Then our Klein geometry can be represented as $(G,H)$.  
We form equivalence classes the same as in (\ref{eq:sen-matrix-equiv}) with the $O$ as in (\ref{eq:fix-subspace-G}).  
Let $H \subset G$ be the subgroup of all matrices of the form
\begin{equation}
\label{eq:fix-subspace-H}
	N = 
    	\left[ 
        	\begin{array}{cc}
         		O 				&	\mathbf{0} \\
            	\mathbf{0}^{T}	&	1
          	\end{array}
    	\right] \mbox{.}
\end{equation}
Then $H$ stabilizes $\mathbf{0} \in \mathcal{Q}$ and we have our Klein geometry as $(G,H)$.  Our manifold is now given by $\mathcal{Q} = G / H$.  
The tangent vectors $q \in \mathcal{T}_{\mathbf{0}}\mathcal{Q}$ have the form
\begin{equation}
\label{eq:subspace-tau}
	q = 
    	\left[
        	\begin{array}{cc}
         		0_{n + m} 		&	Q \mathbf{v} \\
            	\mathbf{0}^{T}	&	0
          	\end{array}
    	\right] \mbox{,}
\end{equation}
and $\exp[q] = I + q$.
Any point $q \in \mathcal{Q}$ can be reached from $\mathbf{0}$ by a path of the form $\exp[q_{1}] \cdots \exp[q_{n}] \mathbf{0}$.  
This is how a path in the subspace would be constructed during an optimization procedure.  

The example is illustrative but not practical.  
No one in their right mind would do this particular optimization problem as above.  
But the methodology presented is generic for Klein geometries associated with a reductive homogeneous space.  
As such it does become a practical method as the Klein geometry becomes more complicated.

\section{\label{sec:lie-second}A Second Look at Lie Groups and Lie Algebras}
The Klein geometries considered above can all be represented as matrices.  
Since they're reductive homogeneous spaces we can remain within the space using the matrix exponential of $\exp: \la{m} \rightarrow G/H$ by restricting the exponential mapping of \la{g} to the subspace \la{m}.  
For a fixed $\omega \in \la{m}$ moving in a specific direction is given by $\exp[t \omega]$, where $t \in \mathbb{R}$.  
It's easier to work with the Lie algebra since it's a vector space.  
For example, one can place a mesh in \la{g} and perform MADS there.

Theorem~\ref{th:lie-path} does not tell us that the exponential map is surjective but we do have that:
\begin{theorem}[Exponential map injectivity radius \cite{book:Stillwell-2008}]
\label{th:exp-diffeo}
The exponential map $\exp: \la{g} \rightarrow G$ is a diffeomorphism in a neighborhood of $0 \in \la{g}$ given by $\| \omega \| < \rho$ for some $\rho > 0$ where $\omega \in \la{g}$.  
The $\rho$ is called the \emph{injectivity radius} of the exponential map.  
Similarly, the matrix logarithm map $\log: G \rightarrow \la{g}$ is a diffeomorphism in a neighborhood of $e \in G$ given by $\| g - e \| < \epsilon$ for some $\epsilon > 0$ where $g \in G$.
\end{theorem}
The $\omega_{i}$ in Theorem~\ref{th:lie-path} would all lie within the injectivity radius of the exponential map.  
Determining the injectivity radius of the exponential map is important since it determines an upper bound on allowable step sizes in any optimization routine.  
This is not always an easy task.  
One may be able to appeal to the group being compact or the space being symmetric to show geodesic completeness in which case Theorem~\ref{th:lie-path} holds with a single $\exp[\omega]$.  
Alternatively, one can develop `tricks' on a case by case basis.

To illustrate the latter consider the special linear group \SLn where $\exp:\sln \rightarrow \SLn$ is not surjective.  
A polar decomposition of a matrix $A \in \SLn$ shows that $A = QS$  where $Q \in \SOn{n}$, $S \in \SPDn$ with $\det(S) = 1$, and both $Q$ and $S$ are unique.  
From Table~\ref{tbl:metrics_exps} we see that every $\omega \in \sln$ can be represented as $\alpha + \sigma$ where $\alpha \in \son{n}$, $\sigma \in \spdn$ with $\mbox{Tr}(\sigma) = 0$, and $\alpha$ and $\sigma$ are also in \sln.

Let \la{m} be the traceless members of \spdn and $\la{h} \equiv \son{n}$.  
Then $\sln = \la{h} \oplus \la{m}$ and $\SLn / \SOn{n}$ is a symmetric space.  
To pick a general member of \SLn first select a $\mu \in \la{m}$ and map this to the base space by $S = \exp[\mu]$.  
This is a surjective map (since \Man is a symmetric space and hence geodesically complete), so we can reach any $S \in \SPDn$, $\det(S) = 1$, from the identity element.  
At $S$ an element $\eta \in \la{h}$ is mapped into \SOn{n} by $Q = \exp[\eta]$.  
This is also a surjective map.  
By forming the product $QS$ we see that every member of \SLn can be reached by this procedure.

Unfortunately this method can not be employed in general.  
The reason it works so well here is that when we decompose a matrix $A \in \SLn$ into the polar form $A = QS$ both $Q$ and $S$ are also in \SLn.  
So we are always working in \SLn and do not have to worry about whether the product $QS$ lies in \SLn or not.

\section{\label{sec:convergence}Algorithm and Convergence}
The most general optimization problem considered here is stated as
\begin{equation}
\label{prob:gen-opt}
    \min_{p \in (G,H)} \ f(p)
    \mbox{ subject to } \mathbf{g}(p) \leq \mathbf{0}
\end{equation}
with $(G,H)$ being a Klein geometry and $\la{g} = \la{h} \oplus \la{m}$.
Any inequality constraints are handled purely within the DSM and will not concern us here.  
They are only allowed if the DSM can handle inequality constraints.
In this section we will initially present a procedure that works with a generic DSM.
This is then specialized to the probabilistic descent method \cite{art:gratton-2015} in Section~\ref{sec:in-group}.
Section~\ref{sec:in-group} also develops a new method of working directly on manifolds.
This new methodology relies on certain properties of the probabilistic descent method.

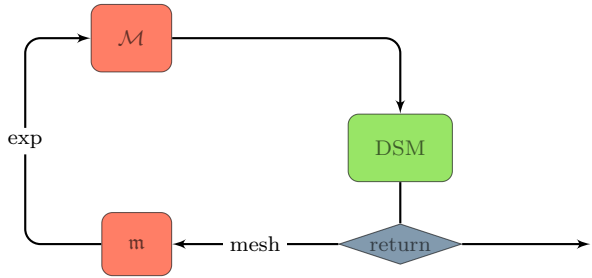
\begin{figure}[htb!]
	\centering
    \begin{tikzpicture}
	    \node [] (dummy1) {}; 
    	\node [special-box, below=1cm of dummy1.center] (node-dsm) {DSM};
    	\node [test-box, below=1cm of node-dsm.center] (dummy2) {return};
        \node [right=2.5cm of dummy2.center] (return) {};
    	\node [dsm-box, left=3cm of dummy2.center] (node-la) {$\la{m}$};
    	\node [dsm-box, left=3cm of dummy1.center] (node-lg) {\Man};
        \draw[-, thick]
        	(node-dsm.south) to (dummy2.north);
		\draw[->, >=latex', thick]
        	(dummy2.west) to node[descr]{mesh} (node-la.east);
        \draw[->, >=latex', thick, rounded corners=2mm]
        	(node-la.west) -- +(-1,0) |- node[descr, pos=.25] {$\exp$} (node-lg.west); 
		\draw[->, >=latex', thick]
        	(dummy2.east) to (return.west);
		\draw[->, >=latex', thick, rounded corners=2mm]
        	(node-lg.east) -| (node-dsm.north);
	\end{tikzpicture}
	\caption{Algorithm with a surjective exponential map.}
    \label{fig:exp-surjective}
\end{figure}

In the simplest case of a surjective exponential map the procedure is done by laying out the mesh in the Lie subalgebra $\la{m} \subset \la{g}$.  
Since the Lie algebra is a vector space this is a trivial mapping.  
The exponential map takes the chosen mesh points and maps them into the reductive homogeneous space $\Man \subset G$.  
These points are then fed back into the DSM.  
(See Figure~\ref{fig:exp-surjective}.)  
Convergence results are the same as whatever are available with the chosen DSM since we have done a global pullback of the objective function (and any inequality constraints) from \Man to \la{m} \cite{misc:dreisigmeyer-2006}.
In this case, ignoring the inequality constraints, (\ref{prob:gen-opt}) becomes
\begin{equation*}
    \min_{\omega \in \mathfrak{m}} \ f \circ \exp (\omega) \mbox{.}
\end{equation*}

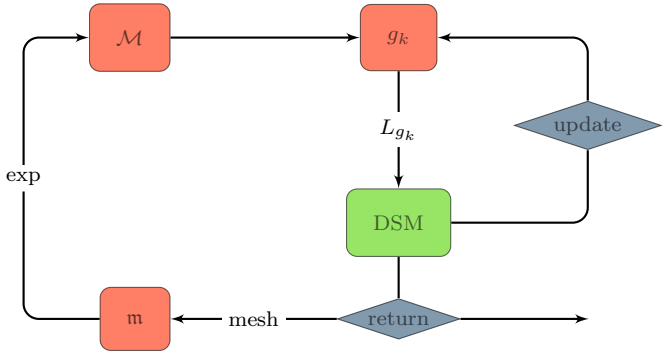
\begin{figure}[hbt!]
	\centering
	\begin{tikzpicture}
	    \node [dsm-box] (dummy1) {$g_{k}$}; 
    	\node [special-box, below=2cm of dummy1.center] (node-dsm) {DSM};
    	\node [test-box, below right=1cm and 2cm of dummy1.center] (node-update) {update};
    	\node [test-box, below=1cm of node-dsm.center] (dummy2) {return};
        \node [right=2.5cm of dummy2.center] (return) {};
    	\node [dsm-box, left=3cm of dummy2.center] (node-la) {$\la{m}$};
    	\node [dsm-box, left=3cm of dummy1.center] (node-lg) {\Man};
        \draw[-, thick]
        	(node-dsm.south) to (dummy2.north);
		\draw[->, >=latex', thick]
        	(dummy2.west) to node[descr]{mesh} (node-la.east);
        \draw[->, >=latex', thick, rounded corners=2mm]
        	(node-la.west) -- +(-1,0) |- node[descr, pos=.25] {$\exp$} (node-lg.west);       
		\draw[->, >=latex', thick]
        	(dummy2.east) to (return.west);
		\draw[->, >=latex', thick]
        	(node-lg.east) to (dummy1.west);
		\draw[->, >=latex', thick]
        	(dummy1.south) to node[descr]{$L_{g_{k}}$} (node-dsm.north);		
    	\draw[-, >=latex', thick, rounded corners=2mm]
			(node-dsm.east) -| (node-update.south);
		\draw[->, >=latex', thick, rounded corners=2mm]
			(node-update.north) |- (dummy1.east);
	\end{tikzpicture}
	\caption{Algorithm without a surjective exponential map.}
    \label{fig:exp-not-surjective}
\end{figure}

With the surjective exponential map the optimization occurs in the fixed tangent (vector) space $\la{m}$.
This is why convergence results from a DSM carry over without modification.
When the exponential map is not surjective more care must be taken.
What we desire is to work in a fixed tangent space after a finite number of steps in the optimization.
If this occurs then convergence results are again the same as the DSM employed.
This situation is shown in Figure~\ref{fig:exp-not-surjective}.  
The only difference is that we may move over \Man, $g_{k - 1}$ and $g_{k}$ may differ, after picking out new elements from \la{m} and mapping them to \Man.
If $g_{k - 1} \neq g_{k}$ then the tangent spaces differ and we essentially start a new optimization with a new mesh.

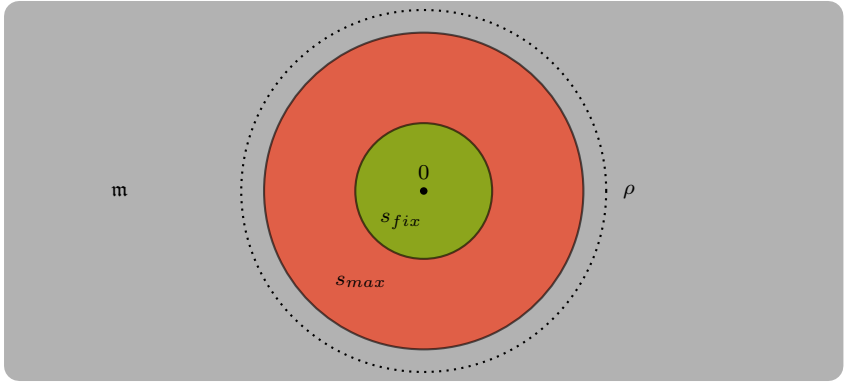
\begin{figure}[hbt!]
	\centering
  	\begin{tikzpicture}
    	[x={(1cm,0cm)}, y={(0cm,1cm)}, z={(0cm,0cm)}, scale=3]
		\begin{scope} [canvas is xy plane at z=0]
			\draw [fill={rgb:red,140;green,0;yellow,26}, opacity=.6, thick] circle (0.7);
        	\draw [fill={rgb:red,34;green,139;yellow,34}, opacity=.6, thick] circle (0.3);
        	\draw [dotted,thick] circle (0.8);
            \node [circle,fill,inner sep=1pt,label=above:$0$] (a) {};
            \node [xshift=-0.3cm, yshift=-0.4cm] (node-fix) { $s_{fix}$};
            \node [below left = 0.4cm and -0.3cm of node-fix] (node-max) { $s_{max}$};
            \node [xshift=-4.0cm] (node-lam) {\la{m}};
            \node [xshift=2.7cm] (node-rho) {$\rho$};
            \node [] (dummy1) {}; 
            \node [above left=2cm and 5cm of dummy1] (dummy2) {}; 
            \node [above right=2cm and 5cm of dummy1] (dummy3) {}; 
            \node [below right=2cm and 5cm of dummy1] (dummy4) {}; 
            \node [below left=2cm and 5cm of dummy1] (dummy5) {}; 
		\end{scope}
        \begin{pgfonlayer}{background}
        	\filldraw [line width=4mm,join=round,black!30]
            	(dummy2.north -| dummy3.east) rectangle 
                (dummy4.south -| dummy5.west);
        \end{pgfonlayer}
	\end{tikzpicture}
	\caption{The area in the vector space \la{m} where the optimization procedure is performed.
    This defines the update procedure in Figure~\ref{fig:exp-not-surjective}.  
    If $0 \leq \|\omega\| < s_{fix}$, $\omega \in \la{m}$, then the tangent space doesn't move.  
    If $ s_{fix} \leq \|\omega\| < s_{max} \leq \rho$, $\rho$ the injectivity radius, then the point $g \in \Man$ is updated.}
    \label{fig:annulus}
\end{figure}

Recall that $\rho$ is the injectivity radius of the exponential map.
For non-surjective exponential maps a maximum step size $0 < s_{max} \leq \rho$ is used.
This guarantees that every point on \Man in a neighborhood of $g_k$ can be reached from $g_k$.
In order to have convergence results we also use a second radius $0 < s_{fix} < s_{max}$ such that we only move to a new point $\exp[\omega] = g_{k} \in \Man$ when the point $\omega \in \la{m}$ lies in the annulus defined by $s_{fix}$ and $s_{max}$.  
(See Figure~\ref{fig:annulus}.)  
Otherwise we set $g_{k} = g_{k - 1}$ and work in the same tangent space.  

Suppose $g_k$ becomes fixed after a finite number of steps in the optimization.
Since we then work in a single vector space where the exponential map is diffeomorphic all the convergence results carry over from the DSM.  
We are now doing a pullback procedure around $g_{k}$ but only locally on \Man.
If $g_k$ never becomes fixed the optimization does not converge.
This is reasonable since for every $g_k$ after a finite number of steps we will move to a new point on \Man that is outside a local neighborhood of $g_k$ determined by $s_{fix}$.

A brief note on how we find the tangent space at each point on \Man.  
Every tangent space looks the same as \la{m} attached at $e \in \Man$.  
If we are at a given point $g \in \Man$ we can work with \TMan{g} by performing the automorphism $g : \Man \rightarrow \Man$.  
That is, we first pick a point $\omega \in \la{m}$, map this to a point $\exp[\omega] \in \Man$ and then move that point by $L_{g}(\exp[\omega]) : \la{m} \rightarrow \Man$.  
This is equivalent to $\exp : \TMan{g} \rightarrow \Man$ and is computationally convenient.

The final optimization routine is given by Procedure~\ref{alg:final_alg}.  
The hardest part of the procedure can be determining what the injectivity radius is if the exponential map is not surjective.  
See the probabilistic descent DSM \cite{art:gratton-2015} in Procedure~\ref{alg:prob-descent-la} for a concrete example.

\begin{procedure}[hbt!]
    \begin{myAlgorithm}[Direct searches on reductive homogeneous spaces]
    \label{alg:final_alg}
    Specify:
    \begin{itemize}
    	\item A direct search method and any required inputs to that method.
        \item A Klein geometry $(G, H)$ with a reductive homogeneous space $\Man = G/H$ and $\la{g} = \la{h \oplus \la{m}}$.
        \item A mapping of the mesh from the vector space \Rn{n} to the vector space \la{m}.
        \item Two radii satisfying $0 < s_{fix} < s_{max} \leq \rho$ if not working with a surjective exponential map.
        Otherwise set $s_{fix} = s_{max} = \infty$.

        \item An initial set $\mu_{0}^{j}$ of mesh points in the Lie algebra \la{m} where $j = 1, \ldots, m$ and $m$ is the number of points required by the DSM.
    \end{itemize}
Let the initial point on the reductive homogeneous space be the identity $g_{0} = e$.  
Finally let $k = 0$.
      \begin{description}[font=\bfseries]
      	\item [Iterate] Set $k \leftarrow k + 1$.
        \item [Move] Find $g_{k}^{j} = L_{g_{k - 1}}(\exp[\mu_{k-1}^{j}])$ for every $j$.
        \item [DSM] Perform the DSM with the $g^{j}_{k}$.
        \item [Return] If the DSM has converged return the optimal point $\widehat{g} \in G/H$ and exit.
        \item [Update] If the current best $\mu_{k}^{*}$ has $\| \mu_{k}^{*} \| > s_{fix}$ update the current point $g_{k} \in G/H$.
        Here $g_{k} = L_{g_{k - 1}}(\exp[\mu_{k}^{*}])$.  
        Otherwise set $g_{k} = g_{k - 1}$.
        \item [Mesh] In \la{m} lay out the new points $\mu_{k}^{j}$ returned by the DSM.
        This is either done in an existing mesh if $g_{k} = g_{k - 1}$ or create a new mesh otherwise.
        \item [Goto] Goto the \textbf{Iterate} step.
      \end{description}
    \end{myAlgorithm}
\end{procedure}

\section{\label{sec:in-group}Working only with the Group}
The Lie algebra occupies a central position in the optimization procedure presented in Procedure~\ref{alg:final_alg}.  
In this way there is no major difference between what was presented here and the method developed in \cite{misc:dreisigmeyer-2006}.  
What is different is the ease that the manifold can be moved over because the symmetry group of the underlying space is known.  
As we stated previously, the reductive condition on a Klein geometry's space $G/H$ is what allows this.

But if the group of symmetries is known and the elements provide a way of moving over the space we may ask: Why do we need the Lie algebra?  
The main reason for the use of a Lie algebra is that as a vector space it is easy to lay out a mesh or pick out a random element.  
If this convenience can be pushed down onto the group itself then the algebra can be dispensed with.  
That creates an alternate way of optimizing in reductive homogeneous spaces.

Let's examine what this means in \Rn{n}.  
Euclidean space is different since it doesn't matter if one looks at the tangent space or the manifold since they both are globally identical.  
In subsection~\ref{subsec:se3} we saw that a point in the Lie algebra and a point on the Lie group can be represented by, respectively,
\begin{equation}
\label{eq:tn-matrix-rep}
	\tau = 
    	\left[
        	\begin{array}{cc}
         		0_{n} 				&	\mathbf{v} \\
            	\mathbf{0}^{T}	&	0
          	\end{array}
    	\right] \mbox{ and }
	\exp[\tau] = 
    	\left[
        	\begin{array}{cc}
         		I_{n} 				&	\mathbf{v} \\
            	\mathbf{0}^{T}	&	1
          	\end{array}
    	\right] \mbox{.}
\end{equation}
For moving in the group, if we work with the Lie algebra addition is used (followed by exponentiation) while the Lie group uses multiplication.  
This gives
\begin{equation}
\label{eq:tn-add-matrix-rep}
	\tau_{1} + \tau_{2} = 
    	\left[
        	\begin{array}{cc}
         		0_{n} 				&	\mathbf{v}_{1} + \mathbf{v}_{2} \\
            	\mathbf{0}^{T}	&	0
          	\end{array}
    	\right] \mbox{ and }
	\exp[\tau_{1} + \tau_{2}] = 
    	\left[
        	\begin{array}{cc}
         		I_{n} 				&	\mathbf{v}_{1} + \mathbf{v}_{2} \\
            	\mathbf{0}^{T}	&	1
          	\end{array}
    	\right] \mbox{.}
\end{equation}
(Of course normally we dispense with all of this and just use $\mathbf{v}_{1} + \mathbf{v}_{2}$ not caring if it is being done in the Lie algebra or the Lie group.)  
The reason the tangent space and the manifold are interchangeable here is that $\exp[\tau_{1} + \tau_{2}] = \exp[\tau_{1}] \exp[\tau_{2}]$ which is due to the fact that $[\tau_{i}, \tau_{j}] = 0$ for any $\tau_{i}$ and $\tau_{j}$.  
A Lie algebra where the Lie bracket is trivial is called \textit{abelian} and has an associated commutative (abelian) connected Lie group.  
What we wish to do here is to work strictly with the $g = \exp[\tau]$ elements in the Lie group and dispense with the elements $\tau$ in the Lie algebra.  
For abelian Lie groups (algebras) such as \Rn{n} and the $n$-dimensional tori this is trivial.  
In fact, $\mathbb{R}$ and the $1$-dimensional torus (the circle) effectively exhaust the abelian Lie groups since any such group is the direct product of a finite number of copies of $\mathbb{R}$ and the $1$-dimensional torus.

For abelian Lie groups it's the case that $\exp[\omega_{1} + \omega_{2}] = \exp[\omega_{1}] \exp[\omega_{2}]$ and so summation in the Lie algebra directly translates to multiplication in the Lie group.  
It follows that picking some point on a mesh in the Lie algebra, which is a sum, easily translates to the Lie group.  
This can be extended in a straightforward way whenever sums in the Lie algebra can be easily translated to products in the Lie group regardless if the group is abelian or not.  
The difficulties arise when this isn't possible.

A probabilistic descent DSM \cite{art:gratton-2015} seems to be particularly well suited to the methodology we have in mind.  
We present an implementation for a reductive homogeneous space in Procedure~\ref{alg:prob-descent-la}.  
This relies on the Lie algebra and the exponential map.

\begin{procedure}[hbt!]
    \begin{myAlgorithm}[Probabilistic descent using the Lie algebra]
    \label{alg:prob-descent-la}
    Specify:
    \begin{itemize}
        \item A Klein geometry $(G, H)$ with a reductive homogeneous space $\Man = G/H$ where $\la{g} = \la{h} \oplus \la{m}$.
    	\item An objective function $f: \Man \rightarrow \mathbb{R}$.
        \item Two radii satisfying $0 < s_{max} < \rho / 2$, $\rho$ the injectivity radius of the exponential map, and $0 < s_{fix} < s_{max}$.
        \item A forcing function $\varrho : \Man \rightarrow \mathbb{R}$ giving sufficient decrease.
        \item A maximum number of function evaluations.
    \end{itemize}
Let the initial point on the reductive homogeneous space be $g_{0}$.  
Set $k = 0$, $g_{0}^{\prime} = g_{0}$, $w_{0} = 0 \in \la{m}$ and $w_{test} = 0 \in \la{m}$.  
Pick a $0 < s_{0} < s_{max}$.
      \begin{description}[font=\bfseries]
      	\item [Choice] Pick a random element $\omega_{k} \in \la{m}$ where $0 < \| \omega_{k} \| < s_{k}$.  
        Set $h^{+} = L_{g_{k}}(\exp[w_{k} + \omega_{k}])$ and $h^{-} = L_{g_{k}}(\exp[w_{k} - \omega_{k}])$.
        \item [Evaluate] Check:
        \begin{description}[align=left,font=\normalfont]
          \item [IF] $f(h^{+}) < f(g_{k}^{\prime}) - \varrho(\|\omega_{k}\|)$ set $g_{k + 1}^{\prime} = h^{+}$ and $w_{test} = w_{k} + \omega_{k}$
          \item [ELIF] $f(h^{-}) < f(g_{k}^{\prime}) - \varrho(\|\omega_{k}\|)$ set $g_{k + 1}^{\prime} = h^{-}$  and $w_{test} = w_{k} - \omega_{k}$
          \item [ELSE] set $g_{k + 1}^{\prime} = g_{k}^{\prime}$ and $w_{test} = 0$
		\end{description}
        \item [Return] If the maximum number of function evaluations has been exceeded return $g_{k + 1}^{\prime}$ and exit.
        \item [Step size] Set $s_{k + 1} = \min(s_{max}, 2 s_{k})$ if the IF or ELIF steps in \textbf{Evaluate} held.  
        Otherwise set $s_{k + 1} = s_{k} / 2$.
        \item [Move] If $\| w_{test} \| > s_{fix}$ set $g_{k + 1} = g_{k + 1}^{\prime}$ and $w_{k + 1} = 0$.  
        Otherwise set $g_{k + 1} = g_{k}$ and $w_{k + 1} = w_{test}$.
        \item [Iterate] Set $k \leftarrow k + 1$.
        \item [Goto] Goto the \textbf{Choice} step.
      \end{description}
    \end{myAlgorithm}
\end{procedure}

In Procedure~\ref{alg:prob-descent-la} it's important to notice that the work is being carried out around the $g_{k}^{\prime} \in \Man$ at each step $k$ even if the current tangent space is being fixed at $g_{k} \in \Man$.  
We fix the working tangent space in order to trivially carry over convergence results.  
It is only moved when there's a danger of leaving the area where the exponential map is a diffeomorphism.  
However, if the procedure is carried out directly on the reductive homogeneous space without reference to \la{m} we can freely move over \Man knowing that at some point the tangent space could be fixed at some point if there is convergence to a solution.  
The only restriction is that our move is never larger than the injectivity radius of the logarithmic map $\log: \Man \rightarrow \la{m}$.  
Then we have the same convergence results as if we worked with \la{m}.  

Assume that we have a black-box matrix generator $\mathbf{B}_{\Man}^{r}$ that can return a random $g \in \Man$ satisfying $\| g - e \| < r$, $0 < r < R$, $R$ a constant.  
Further, $\mathbf{B}_{\Man}^{r}$ must be able to generate all such $g$ that satisfy the norm condition.  
(Note that $R$ can be the injectivity radius of the logarithmic map but this is not a strict requirement.)  
The probabilistic descent algorithm that works solely on \Man is presented in Procedure~\ref{alg:prob-descent}.

\begin{procedure}[hbt!]
    \begin{myAlgorithm}[Probabilistic descent on reductive homogeneous spaces]
    \label{alg:prob-descent}
    Specify:
    \begin{itemize}
        \item A Klein geometry $(G, H)$ with a reductive homogeneous space $\Man = G/H$.
        \item A constant $R > 0$.
        \item A black-box matrix generator $\mathbf{B}_{\Man}^{r}$ that returns a random $g \in \Man$ satisfying $\| g - e \| < r < R$.  
        All $g$ satisfying the norm condition must be returnable by $\mathbf{B}_{\Man}^{r}$.
    	\item An objective function $f: \Man \rightarrow \mathbb{R}$.
        \item A forcing function $\varrho : \Man \rightarrow \mathbb{R}$ giving sufficient decrease.
        \item A maximum number of function evaluations.
    \end{itemize}
Let the initial point on the reductive homogeneous space be the identity $g_{0} = I$.  
Set $k = 0$ and pick a $0 < r_{0} < \rho$.
      \begin{description}[font=\bfseries]
      	\item [Generate] Generate a $M = \mathbf{B}_{\Man}^{r_{k}}$ and set $h^{+} = L_{g_{k}}(M)$ and $h^{-} = L_{g_{k}}(M^{-1})$.
        \item [Evaluate] Check:
        \begin{description}[align=left,font=\normalfont]
          \item [IF] $f(h^{+}) < f(g_{k}) - \varrho(\|M\|)$ set $g_{k + 1} = h^{+}$
          \item [ELIF] $f(h^{-}) < f(g_{k}) - \varrho(\|M^{-1}\|)$ set $g_{k + 1} = h^{-}$
          \item [ELSE] set $g_{k + 1} = g_{k}$
		\end{description}
        \item [Return] If the maximum number of function evaluations has been exceeded return $g_{k + 1}$ and exit.
        \item [Step size] Set $r_{k + 1} = \min(R, 2 r_{k})$ if the IF or ELIF steps in \textbf{Evaluate} held.  Otherwise set $r_{k + 1} = r_{k} / 2$.
        \item [Iterate] Set $k \leftarrow k + 1$.
        \item [Goto] Goto the \textbf{Generate} step.
      \end{description}
    \end{myAlgorithm}
\end{procedure}

We still need to specify the black boxes $\mathbf{B}_{\Man}^{r}$ in Procedure~\ref{alg:prob-descent}.  
The obvious thing to use is the exponential map $\exp : \la{m} \rightarrow \Man$ and pick a random $\omega \in \la{m}$.  
For that we can use the inequality
\begin{equation}
\label{eq:exp-ineq}
\| \exp[\omega] - e \| \leq \| \omega \| e^{\| \omega \|}
\end{equation}
where $\omega \in \la{m}$ satisfies $\| \omega \| < \rho$.  
However, the exponential map is being used here only as a convenience and it is not crucial that this is the $\mathbf{B}_{\Man}^{r}$ in Procedure~\ref{alg:prob-descent}.  
So the algorithm is still independent of \la{m}.  
(As a note, we can use the condition $\| \omega \| < \rho$ since $\exp$ and $\log$ are inverses of each other.  
So working only in the region of \la{m} where $\exp$ is a diffeomorphism is the same as working in the region on \Man where $\log$ is a diffeomorphism.)  
Any other generator is suitable though.  
For example, the unitary, orthogonal and symplectic groups are examined in \cite{art:diaconis1987, art:mezzadri-2007, art:stewart-1980}.  
If we can work with orthogonal groups then the Euclidean group can also be worked with.  
Generating random symmetric positive definite, unipotent, and general and special linear group members are more straightforward prospects.  
The closed form solutions for the exponential maps of the Grassmannian, Stiefel manifold and unipotent group in Table~\ref{tbl:metrics_exps} makes the exponential mapping route more computationally tractable for these manifolds.  
Of course the norm condition still needs to be taken into account.

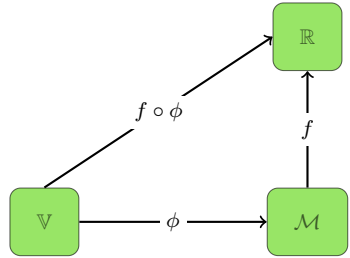
\begin{figure}[htb!]
	\centering
	\begin{tikzpicture}
	    \node [] (dummy1) {}; 
    	\node [special-box, above=1cm of dummy1.center] (node-result) {$\mathbb{R}$};
    	\node [special-box, below=2cm of node-result.center] (node-manifold) {$\Man$};
    	\node [special-box, left=3cm of node-manifold.center] (node-vector) {$\mathbb{V}$};
        \draw[->, thick] 
        	(node-manifold.north) to node[descr]{$f$} (node-result.south);
        \draw[->, thick] 
        	(node-vector.east) to node[descr]{$\phi$} (node-manifold.west);
        \draw[->, thick] 
        	(node-vector.north) to node[descr]{$f \circ \phi$} (node-result.west);
	\end{tikzpicture}
	\caption{Pulling back the function $f : \Man \rightarrow \mathbb{R}$ from the manifold \Man to the vector space $\mathbb{V}$.}
    \label{fig:mesh-distribution}
\end{figure}

One thing we need to consider is the distribution of the random matrices on the manifold \Man.  
Looking at Figure~\ref{fig:mesh-distribution} the pullback procedure is to take the function $f : \Man \rightarrow \mathbb{R}$ and define a new function on a vector space $\mathbb{V}$ by defining $f \circ \phi : \mathbb{V} \rightarrow \mathbb{R}$ where $\phi : \mathbb{V} \rightarrow \Man$.  
As long as $\phi$ is sufficiently smooth this creates no difficulties and the optimization can proceed on $\mathbb{V}$ with all of the convergence results of the employed DSM.  
For a DSM `sufficiently smooth' will typically mean a $\mathcal{C}^{2}$, $\mathcal{C}^{1}$ or Lipschitz function that is at least locally surjective.
For the probabilistic descent method there are conditions on the distribution of the random vectors in the vector space $\mathbb{V}$.  
These are easily enforced, but while the optimization is carried out in $\mathbb{V}$ the final solution actually lies on \Man.  
Given the freedom in choosing $\phi$ it's not obvious how placing restrictions on the vectors chosen in $\mathbb{V}$ carry over to restrictions on the elements of \Man.  
About all that can be stated is that the probability density function of choosing an element must not vanish in any direction from the identity element $e \in \Man$.

\section{\label{sec:applications}Applications}
The procedures developed herein have been applied to two data mining algorithms with positive results.
The first data mining application is concerned with dimensionality reduction of manifold-valued data and involves optimizing over Grassmannians \cite{misc:dreisigmeyer-2017a}.
Also, the objective function is not differentiable and a derivative-free optimization method is the only appropriate procedure.

The second application is an extension of the semi-nonnegative matrix factorization and optimizes over hyper-spheres \cite{misc:dreisigmeyer-2017b}.
This additionally involves the use of optional search steps.
Equality constraints on vector lengths make the restriction to a manifold implicit.
These are naturally enforced by employing the exponential and logarithmic maps from the tangent space to the manifold and vice versa.
That is, at the expense of a more complicated objective function the difficulties of working with equality constraints are removed.
There is also an inequality constraint on the hyper-area defined by the columns of a matrix and it is convenient to treat this as black-box.

Both of the data mining applications have optimization problems that are naturally stated on manifolds.
The exponential and logarithmic maps were used rather than finding points on the manifolds directly.
Modifying the solution method for the modified semi-nonnegative matrix factorization in \cite{misc:dreisigmeyer-2017b} to work solely on the hypersphere \Sn{k - 1}, $k \in \mathbb{N}$, is relativity straightforward.
We show how this can be done presently.

Let $X \in \Rn{n \times m}$ be the data matrix where there are $m$ data points in \Rn{n}, and choose a $k \in \mathbb{N}$.
Typically $k \ll n$.
We're looking for a decomposition $X \approx WH$ where $W \in \Rn{n \times k}$ will have normalized columns, $\|\mathbf{w}_i\| = 1$ for $i = 1,\ldots,k$, and $H \in \Rn{k \times m}$ will have nonnegative entries, denoted by $0 \leq H$.
The measure of model tightness is the maximum geodesic distance
\begin{equation*}
S(W) = \max_{i < j} \arccos( \mathbf{w}_i \cdot \mathbf{w}_j )
\end{equation*}
over \Sn{n-1} between any two columns of $W$.
The optimization problem is
\begin{eqnarray*}
  \min_{W,0 \leq H} \|X - WH\|_{F}^{2} \mbox{ subject to }
  S(W) \leq \epsilon \mbox{ and }
  \| \mathbf{w}_{i} \| = 1 \mbox{.}
\end{eqnarray*}

Take $\widehat{X}$ to be the data set $X$ with columns normalized to unit length and $\bar{\mathbf{x}}$ be the Karcher mean of $\widehat{X}$ \cite{proc:krakowski-2007}.
Find an initial $W_{0} = [ \mathbf{w}_{0, 1} | \cdots | \mathbf{w}_{0, k} ]$ where $S(W_{0}) \leq \epsilon$.
$H_{0}$ is found by solving the nonnegative least squares (NNLS) problem
\begin{equation*}
H_{0} = \mathrm{arg }\min_{0 \leq H} \|X - W_{0}H\|_{F}^{2} 
	\mbox{ where } 
\varepsilon_{0} = \|X - W_{0}H_{0} \|_{F}^{2} \mbox{.}
\end{equation*}

Set
\begin{itemize}
\item $i_{max}$ as the number of iterations to run the algorithm and initialize $i = 0$,
\item $\alpha_{max} = 1$ as the maximum step size,
\item $\alpha_{0} = \alpha_{max}$ as the initial step size,
\item $\theta = 1/2$ as the step size decrease,
\item $\gamma = 2$ as the step size increase, and
\item $\rho(\alpha) = 10^{-3} \alpha^{2}$ as the forcing function.
\end{itemize}
At each step of the algorithm the vectors $\bar{\mathbf{x}}$ and $\mathbf{w}_{i, j}$, $j = 1, \ldots, k$, determine a two-dimensional subspace.
There is a normalized vector $\mathbf{v}_{i, j}$ in this subspace that satisfies $\bar{\mathbf{x}} \cdot \mathbf{v}_{i, j} = 0$, and an angle $0 < \theta_{i, j} < \pi / 2$ such that $\mathbf{w}_{i, j} = \cos(\theta_{i, j}) \bar{\mathbf{x}} + \sin(\theta_{i, j}) \mathbf{v}_{i, j}$.
A \textbf{contraction} of $\mathbf{w}_{i, j}$ is a reduction in the value of $\theta_{i, j}$ by some multiplicative factor $0 < \delta_{i, j} < 1$.
[Below we'll use $\delta = \delta_{i, j} \equiv 0.99$.]
A \textbf{dilation} of $\mathbf{w}_{i, j}$ is an increase in the value of $\theta_{i, j}$ by some multiplicative factor $1 < \gamma_{i, j}$.
Consider an additional normalized random vector $\mathbf{u}_{i, j}$ where $\bar{\mathbf{x}} \cdot \mathbf{u}_{i, j} = \mathbf{u}_{i, j} \cdot \mathbf{v}_{i, j} = 0$, and an angle $0 < \phi_{i, j} < \pi / 2$.
A \textbf{perturbation} of $\mathbf{w}_{i, j}$ is given by
\begin{equation*}
\mathbf{w}_{i, j}^{\pm}
	=
    \cos(\phi_{i, j}) \left[ 
    	\cos(\hat{\theta}_{i, j}) \bar{\mathbf{x}} 
        +
        \sin(\hat{\theta}_{i, j}) \mathbf{v}_{i, j}
	\right]
    \pm
    \sin(\phi_{i, j}) \mathbf{u}_{i, j}
\end{equation*}
where $\hat{\theta}_{i, j}$ will be in some neighborhood of $\theta_{i, j}$.

The optimization routine is as follows:
\begin{description}[font=\bfseries]
	\item [Step 1] Set $i \leftarrow i + 1$.
    \item [Step 2 (Optional Search Steps)]  Try the following two searches.
    \begin{description}[font=\bfseries]
    	\item [Step 2a (optimal solution with contraction)] Set $W^{\prime} = X H^{\dagger}_{i - 1}$, $H^{\dagger}_{i - 1}$ the pseudo-inverse of $H_{i - 1}$, and normalize the columns to unit length.  
    Until $S(W^{\prime}) \leq \epsilon$ and $0 < \theta_{j}^{\prime} < \pi / 2$, iteratively contract the $\mathbf{w}^{\prime}_{j}$ by a factor of $\delta = 0.99$.
    Then find $H^{\prime}$ by solving the NNLS problem.  
    If $\varepsilon_{i-1} - \varepsilon^{\prime} > \rho(\alpha_{i-1})$ set $\varepsilon_{i} \leftarrow \varepsilon^{\prime}$, $W_{i} \leftarrow W^{\prime}$, $H_{i} \leftarrow H^{\prime}$, $\alpha_{i} = \min(\alpha_{max}, \gamma \alpha_{i-1})$ and Goto \textbf{Step 5}.
    	\item [Step 2b (dilation)] Form $\mathbf{w}^{\prime}_{j}$ by dilating $\mathbf{w}_{i - 1, j}$ by a random factor $1 < \gamma^{\prime}_{j} < 1 + \alpha_{i - 1}$ under the restriction $0 < \theta_{j}^{\prime} < \pi / 2$.
If $S(W^{\prime}) \leq \epsilon$ and $\varepsilon_{i-1} - \varepsilon^{\prime} > \rho(\alpha_{i-1})$ set $\varepsilon_{i} \leftarrow \varepsilon^{\prime}$, $W_{i} \leftarrow W^{\prime}$, $H_{i} \leftarrow H^{\prime}$, $\alpha_{i} = \min(\alpha_{max}, \gamma \alpha_{i-1})$ and Goto \textbf{Step 5}.
	\end{description}
	\item [Step 3 (Poll Steps)]  Pick random angles $0 < \theta_{j}^{\prime} < \pi / 2$ where $| \theta_{j}^{\prime} - \theta_{i - 1, j} | < \alpha_{i - 1}$, and $0 < \phi_{j}^{\prime} < \alpha_{i - 1}$.
	\begin{description}[font=\bfseries]
		\item [Step 3a]  Form $\mathbf{w}^{\prime +}_{j}$ by perturbing $\mathbf{w}_{i - 1, j}$.
        If $S(W^{\prime}) \leq \epsilon$ and $0 < \arccos(\mathbf{w}^{\prime +}_{j} \cdot \bar{\mathbf{x}}) < \pi / 2$ then find $H^{\prime}$ and $\varepsilon^{\prime}$ by solving the NNLS problem.
If $\varepsilon_{i-1} - \varepsilon^{\prime} > \rho(\alpha_{i-1})$ set $\varepsilon_{i} \leftarrow \varepsilon^{\prime}$, $W_{i} \leftarrow W^{\prime}$, $H_{i} \leftarrow H^{\prime}$, $\alpha_{i} = \min(\alpha_{max}, \gamma \alpha_{i-1})$ and Goto \textbf{Step 5}.
		\item [Step 3b] Form $\mathbf{w}^{\prime -}_{j}$ by perturbing $\mathbf{w}_{i - 1, j}$.
        If $S(W^{\prime}) \leq \epsilon$ and $0 < \arccos(\mathbf{w}^{\prime -}_{j} \cdot \bar{\mathbf{x}}) < \pi / 2$ then find $H^{\prime}$ and $\varepsilon^{\prime}$ by solving the NNLS problem.
If $\varepsilon_{i-1} - \varepsilon^{\prime} > \rho(\alpha_{i-1})$ set $\varepsilon_{i} \leftarrow \varepsilon^{\prime}$, $W_{i} \leftarrow W^{\prime}$, $H_{i} \leftarrow H^{\prime}$, $\alpha_{i} = \min(\alpha_{max}, \gamma \alpha_{i-1})$ and Goto \textbf{Step 5}.
	\end{description}
	\item [Step 4] Set $\epsilon_{i} \leftarrow \epsilon_{i - 1}$, $W_{i} \leftarrow W_{i-1}$, $H_{i} \leftarrow H_{i-1}$ and $\alpha_{i} = \theta \alpha_{i-1}$.
	\item [Step 5] If $i = i_{max}$ return $W_{i}$ and $H_{i}$, otherwise Goto \textbf{Step 1}.
\end{description}

In the above method the tangent space to \Sn{k - 1} is never used.
Rather, we move directly over \Sn{k - 1} or pick new points on \Sn{k - 1}.
This avoids the need for an exponential (logarithmic) map from (to) the tangent space to (from) the manifold.
The search steps also demonstrate the importance of using the probabilistic descent method in \cite{art:gratton-2015}.
When we let $W = X H^{\dagger}$ and normalize the columns of $W$, the results of mapping these points on \Sn{k - 1} to the tangent space will most likely not lie on any \textit{a priori} existing mesh.

\section{\label{sec:discussion}Conclusions}
If an optimization problem has a feasible set with enough symmetry then it may be possible to find a `law of motion' over the set.
This will have a group structure as well as a notion of distance.
These two requirements, along with an initial feasible point, are all that are needed for performing an optimization routine.
Using the group structure the initial feasible point can be moved over the manifold a specified distance.
This movement allows us to (locally) reconstruct the entire feasible set.
However a difficulty with direct search methods like MADS is the need for a mesh.
The mesh implicitly restricts the group elements and the freedom to reconstruct the feasible set.

When working with nonlinear manifolds, even those with a high degree of symmetry, placing an appropriate mesh on the feasible set can be a daunting challenge.
One way to avoid the problem is to have the mesh be in a tangent space to the manifold.
This does work but the underlying symmetries of the feasible set are lost.
The probabilistic descent method\cite{art:gratton-2015} removes the need for a mesh.
Because of this one is now free to move over the manifold using the `law of motion' with only a restriction on the distance traveled.
Removing the reliance on a mesh allows the symmetries of the feasible set to be fully utilized.



\bigskip
\footnotesize
\textbf{Disclaimer}
Any opinions and conclusions expressed herein are those of the author and do not necessarily represent the views of the U.S. Census Bureau.  
The research in this paper does not use any confidential Census Bureau information.  
This was authored by an employee of the US national government.  
As such, the Government retains a nonexclusive, royalty-free right to publish or reproduce this article, or to allow others to do so, for Government purposes only.

\textbf{Acknowledgements} 
The author would like to thank the editor and reviewers for their helpful comments and suggestions.
\normalsize

\bibliographystyle{spmpsci_unsrt}
\bibliography{jota_rhs}





                       




\begin{landscape}
\renewcommand{\tabcolsep}{3pt}
\begin{table}
\small
\centering
	\begin{tabular}{llllll} 
    	\toprule
			Name & \Man & \TMan{p} & $h(\omega_{1}, \omega_{2})$ & \MyExp{p}{\omega} & Notes \\
        \bottomrule
        \addlinespace
        	Grassmannian & \Gnk{n}{k} & 
            $p^{T} \omega = 0$ &
            $\mbox{Tr}(\omega_{1}^{T} \omega_{2})$ & 
            $\left[ p V \cos(\Theta) + U \sin(\Theta) \right] V^{T}$ & 
            $\omega = U \Theta V^{T}$, \cite{art:edelman-1999} \\
		\addlinespace
        	Stiefel manifold & \Vnk{n}{k} & 
            $p^{T} \omega = - \omega^{T} p$ &
            $\mbox{Tr}\left( \omega_{1}^{T} \left[ I - \frac{1}{2} p p^{T} \right] \omega_{2} \right )$ & 
            $\left[p\ Q \right] 
            \exp\left[
            	\begin{array}{cc}
                p^{T}\omega & -R^{T} \\
                R & 0
            	\end{array}
            \right]
            \left[
            	\begin{array}{c}
                I \\
                0
            	\end{array}
			\right]
            $ & 
            $QR = p_{\perp}p_{\perp}^{T} \omega$, \cite{art:edelman-1999} \\
		\addlinespace
        	General Linear Group & \GLnp & 
            $\GLn{n}$ &
            $\mbox{Tr}(\omega_{1}^{T} \omega_{2})$ & 
            $p \exp(\omega)$ & 
			\pbox{20cm}{$p \in \GLnp$ means\\
            $\det(p) > 0$} \\
		\addlinespace
        	Special Linear Group & \SLn & 
            $Tr(\omega) = 0$ &
            $\mbox{Tr}(\omega_{1}^{T} \omega_{2})$ & 
            $p \exp(\omega)$ & 
            \pbox{20cm}{\MyExp{p}{\omega} is not\\
            surjective} \\
		\addlinespace
        	Special Orthogonal Group & \SOn{n} & 
            $\omega = - \omega^{T}$ &
            $\mbox{Tr}(\omega_{1}^{T} \omega_{2})$ & 
            $p \exp(\omega)$ & 
            See \cite{art:gallier2000} for $\exp(\omega)$ \\
		\addlinespace
        	Special Euclidean Group & \SEn{n} & 
            $\omega = \left[ 
            	\begin{array}{cc}
                A & v \\
                0 & 0
            	\end{array}
            \right]$ &
            $\mbox{Tr}(\omega_{1}^{T} \omega_{2})$ & 
            $p \exp(\omega)$ & 
            \pbox{20cm}{$A = -A^{T}, v \in \mathbb{R}^{n}$ \\ 
            See \cite{art:gallier2000} for $\exp(\omega)$} \\
		\addlinespace
        	\pbox{20cm}{Symmetric Positive Definite\\
            Group} & \SPDn & 
            $\omega = \omega^{T}$ &
            $\mbox{Tr}(p^{-1/2} \omega_{1} p^{-1} \omega_{2} p^{-1/2})$ & 
            $p^{1/2} \exp(p^{-1/2} \omega p^{-1/2})p^{1/2}$ & 
    		\cite{art:pennec-2006} \\
        \addlinespace
        	Unipotent Group & \UPn & 
            $\omega = U , U_{ii} = 0$ &
            $\mbox{Tr}(\omega_{1}^{T} \omega_{2})$ & 
            $\exp(\omega) = \sum_{k = 0}^{n} \omega^{k} / k!$ & 
            $U$ is upper triangular \\
		\addlinespace
		\bottomrule
	\end{tabular}
\caption{The tangent spaces, metrics and exponential maps for some real Lie groups and matrix manifolds.}
\label{tbl:metrics_exps}
\end{table}
\end{landscape}

\end{document}